\documentclass[12pt,a4paper]{article}

\usepackage{amsmath}
\usepackage{amssymb}
\usepackage{latexsym}
\usepackage{srcltx}
 \usepackage{graphics}
\usepackage{graphicx}
\usepackage{epsfig}
\usepackage[ruled,vlined]{algorithm2e}

\newcommand{\re}{{\mathbb R}}

\newcommand{\cA}{{\cal{A}}}
\newcommand{\cO}{{\cal{O}}}

\newcommand{\cS}{{\cal{S}}}

\newcommand{\ba}{{\boldsymbol{a}}}
\newcommand{\bb}{{\boldsymbol{b}}}
\newcommand{\bc}{{\boldsymbol{c}}}
\newcommand{\bd}{{\boldsymbol{d}}}
\newcommand{\be}{{\boldsymbol{e}}}
\newcommand{\bff}{{\boldsymbol{f}}}
\newcommand{\bh}{{\boldsymbol{h}}}
\newcommand{\bn}{{\boldsymbol{n}}}
\newcommand{\bp}{{\boldsymbol{p}}}
\newcommand{\bq}{{\boldsymbol{q}}}

\newcommand{\bu}{{\boldsymbol{u}}}
\newcommand{\bv}{{\boldsymbol{v}}}
\newcommand{\by}{{\boldsymbol{y}}}
\newcommand{\bx}{{\boldsymbol{x}}}
\newcommand{\bz}{{\boldsymbol{z}}}
\newcommand{\bO}{{\boldsymbol{O}}}

\newcommand{\vardot}{\mathord{\,\cdot\,}}

\newtheorem{theorem}{Theorem}
\newtheorem{prop}{Proposition}
\newtheorem{lemma}{Lemma}
\newtheorem{cor}{Corollary}
\newtheorem{remark}{Remark}
\newtheorem{ex}{Example}
\newtheorem{defi}{Definition}

\newtheorem{problem}{Problem}

\date{}

\begin{document}

\author{  
Vladimir Yu. Protasov  
\thanks{University of L'Aquila (Italy)); 
  {e-mail: \tt\small
vladimir.protasov@univaq.it}}}

\title{Eigensets and invariant sets of switching dynamical systems}

\maketitle

\begin{abstract}

We consider reachable sets of switching  systems, which are    
families of linear  ODE ${\dot \bx(t) \, = \, \, A(t) \, \bx(t)}$ with 
a function~$A(\vardot)$ taking values on a given compact set of~$d\times d$ matrices. 
An eigenset is a compact set~$M \ne \{0\}$ that possesses the following property: 
the closure of the set of points reachable  by trajectories~$x(\vardot)$ starting from~$M$ in time~$t$
is equal to 
$e^{\, \alpha t}M, \, t\ge 0$. This concept introduced recently in the literature 
generalizes the notion of an eigenvector of a matrix  to an arbitrary  compact set of matrices. 
We prove the existence of eigensets, analyse their structure and properties, and  
find the corresponding ``eigenvalues''~$\alpha$. The relation of eigensets to the stability of the systems, to their Lyapunov exponents, invariant sets, and 
invariant  norms is established. 
 The question of which compact sets 
can be presented as eigensets of suitable systems is studied. In particular, for~$d=2$, we show that 
every convex~$n$-gone for~$n=3,4,5$, is en eigenset, while for~$n\ge 6$, this is not true.

\bigskip

\smallskip

\begin{flushright}

\end{flushright}

\end{abstract}
\bigskip

\begin{center}
\large{\textbf{1. Introduction}  }
\end{center}
\bigskip

Consider a differential equation on the vector-function~$\bx(t) \in \re^d$: 
\begin{equation}\label{eq.lss}
\left\{
\begin{array}{l}
\dot \bx(t) \ = \ A(t)\bx(t)\, , \qquad t \in [0, +\infty); \\ 
\bx(0) \ = \ \bx_0\, , 
\end{array}
\right. 
\end{equation}
where $A(t)$ is an arbitrary measurable function 
taking values in a given compact set~$\cA$ of real $d\times d$ matrices. 
The set~$\cA$ is called a {\em control set} and the 
  function~$A(\cdot): \, \re_+ \to \cA$ is a {\em switching law}. 
The solutions~$\bx(\vardot)$ called {\em trajectory} 
is an absolutely continuous vector-valued  function~$\bx(\vardot)$ 
satisfying equality~(\ref{eq.lss}) almost everywhere. 
Thus, a switching law and an initial point~$\bx_0$ defines a unique trajectory.

Given a compact set of matrices~$\cA$, we obtain the 
{\em linear switching system}, which  is a family of ODE~(\ref{eq.lss}) with all 
possible switching laws~$A(t)$ taking values on~$\cA$ and all~$\bx_0 \in \re^d$.  
This system is identified with the matrix set~$\cA$. Each matrix~$A\in \cA$ is referred to as 
a {\em regime}. 
If~$\cA$ is a singleton, then we obtain  a linear ODE with constant coefficients. 
Already the case of two-element control set~$\cA$ can offer a great resistance. The equation becomes non-linear since it contains 
the product of two unknowns: the matrix~$A(t)$ and the vector~$\bx(t)$. Arranging  different switches  between  
 regimes from~$\cA$ we obtain trajectories with different growth, different reachable sets, etc.  
Systems with arbitrary control sets~$\cA$ have been studied in an extensive literature, see monographs~\cite{GST12, L03, VdSS} 
and survey~\cite{LA09}.  In particular, the stability of switching systems 
has been a subject of intensive research starting from 1970s due to many applications~\cite{MP86, MP89, O77}, 
see also~\cite{GC, GLP17, ML03} and references therein. 
They  are also analysed in the  optimal control framework as systems with matrix control~\cite{FM12}.  
\begin{remark}\label{r.10} 
{\em Many  works study some  special classes of control sets.  
For example,~\cite{BFRS96, RSS06, V25} consider 
 bilinear systems of the form: 
\begin{equation}\label{eq.lss0}
\left\{
\begin{array}{l}
\dot \bx \ = \ A\bx \ + \ \sum_{i=1}^m u_i B_i\bx \, , \qquad t \in [0, +\infty); \\ 
\bx(0) \ = \ \bx_0\, , 
\end{array}
\right. 
\end{equation}
where for each~$t$, we have~$\bigl(u_1(t), \ldots , u_m(t)) \in U$,  
the set~$U \subset \re^m$ is convex and compact;  
$A, B_1, \ldots, B_m$ are given  $d\times d$ matrices. 
This is a special case of the system~(\ref{eq.lss}) with the convex 
control set~$\cA \, = \, A \, + \, \bigl\{ \, \sum_{i=1}^m u_i B_i : \ \, (u_1, \ldots , u_m) \in U\, \bigr\}$ 
}
\end{remark}

We are interested in the reachable sets and eigensets of linear switching systems.  
The eigensets were introduced in a recent paper of E.Viscovini~\cite{V25} for systems of the form~(\ref{eq.lss0}) under some extra restrictions. An eigenset is  
a compact set~$M\subset \re^d, \, M\ne \{0\}$,  such that for every~$t \ge 0$, 
the closure of the set reachable from~$M$ in time~$t$ coincides with~$\, e^{\, \alpha t} M$, where~$\alpha \in \re$ is a constant. For a system with  a convex control set~$\cA$ the word 
 ``closure'' can be omitted. Thus, an eigenset reproduces itself and multiplies by~$e^{\alpha t}$ 
 at every moment~$t$. 
 For a one-element system~$\cA = \{A\}$, each eigenvector with a real eigenvalue 
 is an eigenset. 
 In some sense, eigensets generalize the concept of eigenvector 
 to arbitrary compact set of matrices.    
We prove the existence of eigensets for arbitrary system with compact 
control set and find the eigenvalues. In particular, 
if the system is irreducible, then the eigenvalue~$\alpha$ 
is the same for all eigensets and is equal to the Lyapunov exponent. 
 Then we find the general structure of eigensets and establish a relation 
 to invariant (Barabanov) norms. A special attention is paid to the geometry of convex eigensets 
  (Section~6). It is shown that every 
 convex eigenset~$M$ contains the origin. Moreover, the origin can either be  an interior point of~$M$ or a conical point of~$\partial M$. The case when~$M$ is a polytope ({\em eigenpolytope}) 
 is analysed. 
 In Section~7 we raise the issue of classifying convex sets that can be realized
 as eigensets of suitable systems. For some sets (symmetric about the origin, arbitrary simplices and parallelepipeds, 
 ellipsoids for~$d\ge 3$) such systems always exist. On the other hand, there is no system in~$\re^2$
 whose eigenset is a disc centered not in the origin.  Finally, Section~8 focuses on 
planar systems ($d=2$).   The main result states that for every convex flat $n$-gon, where~$n\le 5$, 
there exists a system whose eigenset is equal (isometric) to that~$n$-gon, while for~$n=6$, this is not true.

We use the standard notation: ${\rm co}\,(X), \, {\rm int}\, X, \, \partial X$ is, respectively, 
the convex hull, the interior, and the  boundary of the set~$X$, $(\bx, \by)$ 
is the standard scalar product is~$\re^d$, 
 $\rho(A)$ is the spectral radius of the matrix, i.e., the largest modulus of its eigenvalues. 
The default norm in~$\re^d$ is Euclidean. The basis in~$\re^d$ 
is assumed to be fixed, and thus, we identify matrices with the corresponding linear operators; 
$I$  denotes the identity matrix. 
We  denote vectors  by bold letters and scalars by standard letters.    
 
\bigskip

 \begin{center}
\large{\textbf{2. The statement of the problem}  }
\end{center}
\bigskip

For linear switching system~(\ref{eq.lss}) with a convex control set~$\cA$ 
and for an arbitrary compact set~$M \subset \re^d$, the notation~$M_t$ 
denotes
the set of points reachable from the set~$M$ in time~$t$. 
Thus, for every~$\by \in \re^d$, we have~$\by \in M_t$ if and only if there exists a trajectory 
for which~$\bx(0) \in M, \bx(t) \, = \, \by$. We call~$M_t$ the {\em reachable set in time~$t$}.
In particular,~$M_0=M$. Since the control set~$\cA$ is convex, it follows that  all reachable sets are compact (see Remark~\ref{r.100} in Appendix). Actually,   many of the results  obtained in this paper do not actually
rely on the closeness of~$M_t$, the existence of converging sequences of trajectories suffices. 

If~$\cA$ is not convex, then we denote by the reachable set~$M_t$ the closure 
of points reachable from~$M$ in time~$t$. It coincides with the reachable set 
for the system~${\rm co}(\cA)$, since any of its trajectories can be well approximated by trajectories 
of the system~$\cA$, see Remark~\ref{r.100}. By default we assume that~$\cA$ is convex and does not contain zero.

\begin{defi}\label{d.10}
A compact set~$M \subset \re^d$ different from the one-point set~$\{\bO\}$ (the origin)
is called an eigenset for the system~(\ref{eq.lss}) if there exists a constant~$\alpha \in \re$ 
such that~$M_t =  e^{\, \alpha t}M$ for every~$t \in [0, +\infty)$. 
In this case~$\alpha$ is called an eigenvalue. 
\end{defi}
Our aims are to prove the existence of eigensets,  establish their structure and 
topological properties for arbitrary system~(\ref{eq.lss}), find the eigenvalues~$\alpha$. 
Then  we attack the {\em realizability problem}:  which compact sets~$M$  are eigensets of suitable 
systems~(\ref{eq.lss})? The results are especially interesting if~$M$ 
is a convex body or polyhedron. Even for flat polygons the answer is nontrivial and, in a sense, surprising (Section~8).  

An eigenset is never unique: if~$M$ is an eigenset, then so are all~$\lambda M, \, \lambda \ne 0$. 
It may not be unique up to multiplication  by a constant. 
The closure of the union of a bounded (contained in some ball)
family of eigensets is also an eigenset. Note that a convex hull of eigensets may not be an eigenset. 
In general,  the eigensets can 
have a nontrivial structure. They are not necessarily convex and not necessarily connected.
They can have a fractional Hausdorff dimension. Already for second-order systems, different cases are possible.  
For a wide class of systems in~$\re^2$, namely, 
for systems without dominant regime,  the eigenset is always unique and convex (Theorem~\ref{th.80}). 
On the other hand, there exists a system in~$\re^2$ that has no convex eigensets but has infinitely many 
non-convex ones (Example~\ref{ex.20}). 

A weaker version of eigensets are {\em invariant sets},  which are defined in the next section.

\newpage

 \begin{center}
\large{\textbf{3. Eigenvalues vs Lyapunov exponents}  }
\end{center}
\bigskip

The  {\em Lyapunov exponent} of the system~(\ref{eq.lss}) is the maximal exponent of growth  
of its trajectories: $\sigma(\cA)\, = \, \inf\, 
\{\beta \in \re \, : \ \|\bx(t)\| \, \le \, C\, e^{\, \beta t} \ 
\}$, where $C$ 
is a constant depending on the trajectory and  the infimum is computed over all trajectories~$\bx(\vardot)$
generated by all possible switching laws.  For one operator~$A$, the Lyapunov exponent~$\sigma(A)$
is equal to the maximal spectral abscissa, i.e., the maximal real part of eigenvalues. 
A {\em shift} of the system~$\cA$ by a number~$\beta$ is the system~$\cA - \beta I \, = \, \{A-\beta I: \ A \in \cA\}$, where~$I$ is the identity matrix. 
If~$\bx(t)$ is a trajectory of~$\cA$, then~$e^{-\beta t}\bx(t)$ is a trajectory for~$A-\beta I$. 
Therefore, all shifts of the system have the same eigensets. On the other hand, 
taking~$\beta = \sigma(\cA)$, we obtain a {\em normalized} system, for which the Lyapunov 
exponent is equal to zero. Thus, it is sufficient to study the 
eigensets for normalized systems only. 

 \begin{defi}\label{d.20}
System~(\ref{eq.lss}) is called {\em normalized} if~$\sigma(\cA) = 0$. 
The system is called {\em irreducible} if the matrices of~$\cA$ 
do not share a common invariant nontrivial linear subspace. 
\end{defi}
We are going to show that all eigenvalues of an irreducible system~$\cA$ 
are equal to~$\sigma(\cA)$. This will be done by applying  invariant sets 
and invariant norms. We begin with definitions. 
\begin{defi}\label{d.15}
A compact set~$K \subset \re^d$ different from~$\{\bO\}$ 
is {\em invariant}  for the system~(\ref{eq.lss}) if~$K_t \subset   e^{\, \sigma t}K$ for every~$t \in [0, +\infty)$. 
\end{defi}
If~$\cA$ is normalized, then we have~$K_t \subset K$ for all~$t$. 
Thus, every trajectory starting in the invariant set~$K$ never leaves~$K$.
\begin{defi}\label{d.17}
Let~$\cA$ be a system and $K\subset \re^d$ be an arbitrary set. 
A point~$\by$ is  reachable 
from~$K$ in arbitrary time, if for every~$T \ge 0$, there is a
trajectory of~$\cA$ such that~$\bx(0) \in K$ and~$\bx(T) = \by$. 
\end{defi}
We see that a compact set~$M\ne \{\bO\}$ is an eigensest of a normalized system~$\cA$ if and only if~$M_t \subset M$ for all~$t \ge 0$ and all points of~$M$ are reachable from~$M$ in arbitrary time. 
Each invariant set satisfies that first condition by may not satisfy the second one. Thus, {\em an eigenset is 
 invariant, but not vice versa.}

Among all invariant sets, one class is especially important. 
Those are balls of invariant norms. For a normalized system~(\ref{eq.lss}), 
a norm~$f(\vardot)$ in $\re^d$ is called {\em invariant} or 
{\em Barabanov} if a) the function~$f(\bx(t))$ is non-increasing in~$t$ along every trajectory~$\bx(\cdot)$; b)
for an arbitrary point~$\bx_0 \in \re^d$, there is an {\em extremal trajectory} $\bar \bx(\vardot)$n  starting 
at~$\bx_0$, for which~$f(\bar \bx(t)) \,  = \, {\rm const}, \, t \in [0, + \infty)$.  
The theorem of Barabanov~\cite{B88} asserts that 
{\em every irreducible system with a convex control set possesses an invariant norm}. 

An invariant norm of an arbitrary system~$\cA$ is defined as an invariant norm of the corresponding 
normalized system~$\cA - \sigma I$. Thus, 
for an arbitrary trajectory, the function~$e^{-\sigma t}f(\bx(t))$  is non-increasing and, for an extremal 
trajectory, it is constant.   
\begin{remark}\label{r.20} 
{\em Let~$B=\{\bx \in \re^d: \ f(\bx)\le 1\}$ be the unit ball of the invariant norm. The invariant ball is an invariant set but not necessarily an eigenset of the system. 
}
\end{remark}

Now we are ready to find eigenvalues of an irreducible system.

\begin{theorem}\label{th.10}
For every eigenset of an irreducible system, we have~${\alpha = \sigma(\cA)}$. 
\end{theorem}
 {\tt Proof}. Consider an arbitrary invariant norm~$f$ of the system~$\cA$ and denote by~$\bx_0$ 
 the point of an eigenset~$M$, where the maximal value of the norm~$f$ is attained. 
For an arbitrary extremal trajectory~$\bar \bx (\vardot)$ starting at~$\bar \bx(0) = \bx_0$, we have  $f(\bar \bx(t)) \, = \, e^{\sigma  t}f(\bx_0)$. On the other hand, 
$\bar \bx(t) \in e^{\alpha t}M$, therefore~$f(\bar \bx(t)) \, \le \, e^{\alpha t} \max_{\by \in M} f(\by)\, = \, e^{\alpha t} f(\bx_0)$.  
Thus, $e^{\alpha t} \ge e^{\sigma t}$, hence, $\alpha \ge \sigma$. 
Furthermore, for every~$T>0$, we have~$e^{\, \alpha T}\bx_0 \in e^{\, \alpha T}M = M_T$, consequently, 
there exists a trajectory~$\by(\vardot)$ for which 
$\by(0) \in M$ and~$\by(T)\, = \, e^{\, \alpha T}\bx_0$. On the other hand,~$f(\by(T)) \,  \le \, e^{\, \sigma T} f(\by(0))\, \le \, 
e^{\, \sigma T} f(\bx_0)$. 
Thus, $e^{\alpha T}\, \le \,  e^{\sigma T}$, which implies that~$\alpha 
\le \sigma$.

{\hfill $\Box$}
\medskip

\begin{cor}\label{c.10}
If a system is normalized and irreducible, then for each of its eigensets~$M$, 
we have~${\alpha = 0}$. Therefore,~$M_t = M, \ t\in [0, +\infty)$.  
\end{cor}
We see that an eigenset of a normalized system reproduces itself at every moment of time. 
\smallskip

Note that we have not yet proven that an irreducible system possesses an  eigenset. 
Theorem~\ref{th.10} does not imply the existence, it asserts that if there are eigensets then 
they all have the same eigenvalue equal to the Lyapunov exponent. The existence will be established in the next section. We begin with Theorem~\ref{th.20} according to which every system~$\cA$ 
 does possess 
 at least one eigenset with the eigenvalue~$\alpha = \sigma(\cA)$. If~$\cA$ is reducible, 
 that it  can also have 
 other eigenvalues; all of them  are than~$\sigma(\cA)$.  Then we stablish the general structure of eigensets (Theorem~\ref{th.30}).

\bigskip 	 

\begin{center}
\large{\textbf{4. The existence and properties of eigensets}}   
\end{center}

\bigskip

 The  {\em fundamental matrix}~$\Pi(t)$
of a switching law~$A(\vardot)$ 
is a solution of the matrix ODE
$$
 \frac{d}{dt}\, \Pi(t)\ = \ A(t)\Pi(t)\, , \qquad  \Pi(0) = I\, .
$$ 
The corresponding trajectory of the system~(\ref{eq.lss}) is~$\bx (t) = \Pi(t)\bx(0)$. 
Each matrix~$\tilde A\in \cA$ defines the  {\em stationary} switching law
~$A(t)\, \equiv \, \tilde A, \, t \in [0, +\infty)$. For the sake of simplicity, we 
denote it by the same symbol as the matrix and call it {\em the switching law~$\tilde A$.}  
We say that a trajectory~$\bx(\vardot)$ is 
 obtained {\em by applying the switching law~$A(\vardot)$ to the point~$\bx_0$. }
In particular, if we  apply the regime~$\tilde A$ to~$\bx_0$, then we obtain the trajectory
~$\bx(t) \, = \, e^{\, t \tilde A }\bx_0$.

 \begin{theorem}\label{th.20}
Every linear switching system in~$\re^d$ possesses an eigenvalue~$\sigma(\cA)$
and can also have 
 at most~$d-1$ other eigenvalues, which are smaller than~$\sigma(\cA)$. 
\end{theorem}
By Definition~\ref{d.10}, each eigenvalue of a system~$\cA$  corresponds to at least one eigenset.  
Theorem~\ref{th.20} contains two assertions: 1) a system~$\cA$ always has an eigenset with the 
eigenvalue~$\sigma(\cA)$; 2) there are at most~$d$ eigensets with different eigenvalues; all of them 
are smeller than~$\sigma(\cA)$.  

\begin{cor}\label{c.15}
An irreducible system has a unique eigenvalue, which is the Lyapunov exponent. 
\end{cor}

Proof of Theorem~\ref{th.20} is in Appendix. 

\medskip 
 
Now we address the general structure of eigensets.   We are going to show that every eigenset 
 is a union of  ``elementary'' eigensets called {\em ivies.}  Each ivy 
grows from an infinite in both directions trajectory located on the invariant sphere~$\partial B$
(sphere of the invariant norm). 
 
 We call  {\em stem} an infinite in both directions trajectory, 
 i.e., an absolutely continuous map~$\bx: \re \to \re^d$ 
such that for almost all $t\in \re$, there is  $A(t) \in \cA$ for which~$\dot \bx(t)\, = \, A(t)\bx(t)$. 
The stem will be associated with the image of the map~$\{\bx(t): \, t\in \re\}$. 
The {\em ivy} is the closure of all points reachable from the stem.
For a normalized system, we call a stem  {\em extremal}
 if it is entirely located on a sphere of the invariant norm, i.e., if  $f(\bx(t)) \, \equiv \, {\rm const}, \, t\in \re$. 
The corresponding ivy is also called extremal. Let us remark that the points of an extremal 
ivy can lie not only on the sphere but also  inside it.

 \bigskip

 \begin{theorem}\label{th.30}
Every extremal ivy  is an eigenset. 
Conversely, every eigenset of an irreducible normalized system is a union of some set of extremal ivies.  
\end{theorem} 

\begin{lemma}\label{l.10}
 Let~$M$ be an invariant set of a normalized system~$\cA$, 
 $\bx(\vardot)$ is arbitrary trajectory of~$\cA$.
If the point~$\bx(0)$ is reachable from~$M$ in arbitrary time, then 
so is every point~${\bx(t), \, t \ge 0}$. 
 \end{lemma}
{\tt Proof}. For every~$\tau \le t$, the point~$\bx(t)$ is reachable in time~$\tau$ from
the point~$\bx(t-\tau) \in M$. If~$\tau > t$, then we take~$\by \in M$ for which~$\bx(0) \in M_{\tau - t}(\by)$
(it exists by the assumption); then~$\bx(t)$ is reachable in time~$\tau$ from~$\by$. 
 
{\hfill $\Box$} 
\medskip  

 {\tt Proof of Theorem~\ref{th.30}}. Sufficiency follows from Lemma~\ref{l.10}. 
 Indeed, each point~$\bx$ of the steam  is reachable in arbitrary time, hence, so are all 
 points reached from~$\bx$, i.e., all points form the ivy. 
 To prove the necessity, we  consider an invariant norm~$f$ of 
 an irreducible normalized system~$\cA$ and,  
for an arbitrary point~$\bx \in M$,  consider the value
 $$
 r_{\bx} \ = \ \sup\, \bigl\{f(\bu):\  \bu \in M, \ 
 \bx \ \mbox{is reachable from} \  \bu \bigr\}\, .
 $$
It is clear that~$r_{\bx} \ge f(\bx)$. If~$r_{\bx}=f(\bx)$, then take arbitrary~$h>0$
and consider the sequence of points~$\bx_k \in M$ such that
 $\bx_0 = \bx$ and~$\bx_{k}$ is reachable from~$\bx_{k+1}$ in time~$h$. 
Since the sequence~$f(\bx_k)$ is non-decreasing and its supremum does not exceed~$r_{\bx} = f(\bx_0)$, 
we see that it is identically equal to~$r_{\bx}$. Thus, all the points~$\bx_k$ lie 
on the sphere of invariant norm of radius~$r_{\bx}$. Then all the corresponding finite trajectories 
that consecutively connect those points also lie on the same sphere.  
We add the extremal trajectory~$\bar \bx(\vardot)$ starting at~$\bar \bx(0) = \bx$
and obtain infinite in both directions trajectory on the unit sphere. This is en extremal stem
 passing through~$\bx$.

Assume now that~$r_{\bx} > f(\bx)$. Choose a sequence of points~$\by_k \in M$
such that the point~$\bx$ is reachable from each point~$\by_k$, i.e., 
 $\Pi_k \by_k  = \bx$ for some fundamental matrices~$\Pi_k$,  and
  $\lim_{k\to \infty} f(\by_k)\, = \, r_{\bx}$. 
Passing to a subsequence, we assume that~$\by_k$ converges to some point~$\by\in M$.
Clearly~$f(\by) = r_{\bx}$.  
Let~$\bz_k = \Pi_k \by$. We have 
 $$
 f(\bz_k - \bx)\ = \ 
 f\Bigl(\Pi_k\by - \Pi_k\by_k\Bigr)\, = \, f\Bigl(\Pi_k(\by -\by_k)\Bigr)\, \le \, 
 f(\by -\by_k)\, , 
 $$
 since the operator~$\Pi_k$ does not increase the invariant norm. 
Thus, $f(\bz_k - \bx) \, \le \, f(\by -\by_k)$ for all~$k$, 
Consequently,  $\bz_k \to \bx, \, k\to \infty$. 
In this case~$\bx$ is a limit point of trajectories going from the point~$\by$. 
It remains to show that~$\by$ belongs to some extremal stem. 
The supremum of norms~$f(\bu)$ over all points~$\bu\in M$ from which one can reach~$\by$ 
is equal to~$r_{\by}$. On the other hand, $r_{\by} = r_{\bx}$, 
 otherwise there would exist a point~$\bu$, from which every point~$\bz_k$ is reachable 
 and~$f(\bu) > r_{\bx}$. Thus, $r_{\by} = r_{\bx}$. Then by what it is shown above,~$\by$ 
 lies on the extremal stem.

{\hfill $\Box$}
\medskip

  \bigskip 	 

\begin{center}
\large{\textbf{5. Invariant sets  of generic systems \\ contain zero and are arcwise connected} } 
\end{center}

\bigskip  
 
We begin the analysis of the topology of eigensets with 
proving that they are arcwise connected, apart from one special case. 
The  only exception  is given by isotropic systems, which have proportional
(in suitable basis) orthogonal Lie group. In this basis,   
 the only invariant norm is Euclidean, 
 and every extremal ivy (Theorem~\ref{th.30}) lies on a Euclidean sphere. 
\begin{defi}\label{d.30}
A linear switching system~$\cA$ is called {\em isotropic} if there is a basis in~$\re^d$
in which every matrix from~$\cA$ gets the form~$A \, = \, s I + X$, 
where~$X$ is an skew-symmetric matrix and the number~$s$ is the same for all~$A \in \cA$. 
Such bases are also called isotropic. 
\end{defi} 
Every isotropic system in the corresponding basis 
has the Lie algebra proportional to the orthogonal Lie algebra,   and vice versa. 
This terminology originated with  ``isotropic matrix'', a diagonalizable matrix with 
equal by modulus eigenvalues,  often used in the literature. 
If the system is isotropic, then the exponents of all its matrices are isotropic.  An arbitrary skew-symmetric matrix~$X$ have all its exponents~$e^{tX}$ being 
 orthogonal matrices. Therefore, $e^{tA} = e^{t s}U(t)$, where~$U(t)$ is orthogonal. 
All trajectories of this system satisfy $\| \bx(t)\| \, = \,  e^{\, st}\|\bx(0)\|$, 
and therefore,~$s=\sigma(\cA)$. In particular, for a normalized system, every trajectory 
lies on a Euclidean sphere. This implies the following 
  \begin{prop}\label{p.35}
For an arbitrary isotropic system~$\cA$, we have~$s = \sigma(\cA)$. 
If~$\sigma(\cA) = 0$, then in a suitable basis, all trajectories of~$\cA$ are located on Euclidean 
spheres centered at the origin. An intersection of each eigenset/invariant set of~$\cA$ with an arbitrary Euclidean sphere centered at the origin is also an eigenset (respectively, invariant set), unless it is empty. 
\end{prop}
Applying Theorem~\ref{th.20} we obtain  
\begin{cor}\label{c.20}
For an isotropic system, all invariant sets in the isotropic basis are closures 
of the unions~$\bigcup_{r \ge  0} r M^{(r)}$, where each~$M^{(r)}$ is either empty or 
an invariant set on the unit sphere with center at the origin. The same is true for eigenset, 
in this case each~$M^{(r)}$ is either empty or 
an eigenset. The same is true for invariant sets. 
\end{cor}
Thus, an isotropic system possesses an eigenset~$M$ on the surface of a Euclidean sphere (in the 
isotropic basis).  
The sets~$\lambda M, \, \lambda \in \re$,  can compose compact disconnected sets as well as sets with a fractional 
Hausdorff dimension. It turns out, however, that this is possible only for isotropic systems.   
\begin{theorem}\label{th.50}
Every invariant set of a non-isotropic system contains zero. 
\end{theorem}
The proof is based on properties of multiplicative matrix semigroups with constant spectral radius. 
They have been intensively studied, see~\cite{CMMS23, OR97, P13} and references therein. We will rely on the following theorem proved in~\cite{PV17}: 
 \medskip 

\noindent \textbf{Theorem A}. {\em Let~$\cS$ be an irreducible multiplicative semigroup
of~$d\times d$ matrices. 
Then the following folds:  all nonsingular matrices of~$\cS$ 
have the same spectral radius if and only if  there is a basis in~$\re^d$ in which they are all orthogonal.}
\medskip 

Since all the matrices~$e^{\,t\,A}, \, A\in \cA$,  are nonsingular, 
they generate a semigroup  to which 
Theorem~A can be applied. 
  \begin{prop}\label{p.40}
All fundamental matrices~$\Pi(t)$ of an irreducible systerm~$\cA$ have the same spectral radius if and only if  there is a basis in which  all matrices of~$\cA$ are skew-symmetric. In this case the spectral radius is equal 
to~$1$. 
\end{prop}
The sufficiency in Proposition~\ref{p.40} is simple. If, for all~$A\in \cA$, it holds that~$A^T = - A$, 
then $(A\bx , \bx) = 0$ for all~$\bx \in \re^d$, therefore, every trajectory~$\bx(\vardot)$
has the velocity at each point orthogonal to the radius-vector. Therefore, $\|\bx(t)\|\, \equiv \, {\rm const}$, 
and hence, all the fundamental matrices~$\Pi(t)$ are orthogonal. The necessity is less obvious, since all the matrices from~$\cA$ must become skew-symmetric simultaneously in a common basis. 
\smallskip 

 {\tt Proof of Proposition~\ref{p.40}}.  All fundamental matrices~$\Pi(t)$ 
 are nonsingular and form a multiplicative semigroup. By Theorem~A, 
 there exists a basis in which all of them are orthogonal. 
Hence, in this basis, we have  $\|\bx(t)\| \, \equiv \, {\rm const}$ along every trajectory~$\bx(\vardot)$. Consequently, for every switching law, it holds 
that~$\bigl(A(t)\bx , \bx \bigr)\, = \, \bigl(\dot \bx , \bx \bigr)\, = \, \frac12 \, \frac{d}{dt}\, \bigl(\bx, \bx \bigr)\, = \, 0$. Thus, $\bigl(A\bx, \bx \bigr) = 0$ 
for each point of the trajectory~$\bx$. Since the initial point is arbitrary, 
we see that the equality~$(A\bx, \bx) = 0$ is true for all~$\bx \in \re^d$, 
therefore,~$A^T = -A$. 
Conversely, if~$A^T = -A$ for all~$A\in \cA$, 
then the norm 
 $\|\bx(t)\|$ is constant along each trajectory. 
 This yields that all the matrices~$\Pi(t)$ 
 are orthogonal and have spectral radius equal to one. 
 
{\hfill $\Box$}
\medskip  

Now we begin the proof of the fundamental theorem. 

\medskip

\noindent {\tt Proof of Theorem~\ref{th.50}}. 
We assume that the system~$\cA$ is normalized  and that~$B$ is a unit ball
of Barabanov's norm. 
Since~$\Pi(t) B \subset B$ and therefore,~$\Pi^k(t) B \subset B$ for all~$k$,
it follows  that~$\rho (\Pi(t)) \le 1$.  
If there is~$t$ for which~$\rho (\Pi(t)) < 1$,  then for the switching law 
with period~$\Pi(t)$, we have~$\Pi(kt) = \Pi^k(t) \to 0$
as~$k \to \infty$ and hence,~$\Pi(kt)M \to 0$. 
Thus,  every invariant set~$M$ contains zero. If, otherwise, 
 all the matrices~$\Pi(t), \, t\ge 0$, have the spectral radius~$1$, 
then, invoking  Proposition~\ref{p.40}, we conclude that all~$A \in \cA$ 
are antisymmetric in a common basis, hence, the matrices~$\Pi(t)$ are orthogonal. 
In this case every trajectory of the system lies on a Euclidean sphere, which completes the proof.

{\hfill $\Box$}

As a corollary, we conclude the connectedness of an arbitrary eigenset of a non-isotropic system.

  \begin{theorem}\label{th.52}
All invariant sets of a non-isotropic system are arcwise connected. 
\end{theorem}
{\tt Proof}. If a normalized system is not isotropic, then 
there is a fundamental matrix~$\Pi(T)$ with the spectral radius 
not equal to one, and hence, smaller than one. The corresponding periodic switching
law with period~$T$ moves an arbitrary point~$\bx_0 \in M$
along a trajectory that converges to zero. Therefore, 
$\bx_0$ is connected with the origin  by this trajectory. 
 
{\hfill $\Box$} 


 \bigskip 	 

\begin{center}
\large{\textbf{6. Convex eigensets and polytopes} }
\end{center}

\bigskip

Not all eigensets are convex; there are systems 
with no convex  eigensets at all (see examples in  Section~9). 
It is a challenging problem to characterise 
systems that have at least one convex eigenset. Note that 
for invariant sets, this issue is simple and  the answer is always affirmative: 
the unit ball of Barabanov's norm is a convex invariant set. 
Convex eigensets 
possess special geometrical properties. To formulate the main results we need some further notation. 

A convex body is a convex compact set with a nonempty interior. 
An outer normal vector to a convex body~$G$ at a point~$\ba \in \partial G$
is a nonzero vector~$\bn$ such that~$\max_{\bx \in G} (\bn, \bx - \ba) \, = \, 0$. 
All normal vectors span a {\em normal cone}. A point~$\ba$ is called {\em conical} if there exists an outer normal vector~$\bn$ and $c> 0$ such that~$(\bn, \bx - \ba) \, \ge  \, c\, \|\bx - \ba\|$
for all~$\bx \in G$. Equivalently,  the {\em tangent cone} generated  by the vectors~$\{\bx - \ba\, : \ \bx \in G\}$ is pointed, i.e., 
does not contain a straight line. The 
normal cone at a conical point  has dimension~$d$. Every conical point of a convex body is an extreme point, 
but not vice versa. 

As usual, a face of~$G$ is a convex subset~$F \subset G$ such that, for any
line segment $[a, b] \subset  G$  such that $F \cap  (a, b) \ne \emptyset$, we have 
$[a, b] \subset F$. The face of dimension~$d$ is $G$ itself, faces of dimensions~$d-1$ and~$1$ are, respectively, 
facets and extreme points. A {\em carrier} of a point~$\bx \in G$ is the smallest by inclusion 
 face containing~$\bx$. The dimension of the carrier is denoted by~$\delta(\bx)$.

 \bigskip 

\begin{center}
\textbf{6.1. The main results} 
\end{center}

\bigskip 
  
If a system is irreducible, then every invariant set must be full-dimensional, i.e., 
is not contained in a proper linear subspace. Hence, every convex invariant set is a convex body. 
We begin with a necessary condition for invariant convex bodies.  

 \begin{theorem}\label{th.54}
Every convex invariant set~$M$ of an irreducible system contains the origin. In this case, the origin 
can be  either an interior point or a conical point of~$M$.  
\end{theorem}
The first part of Theorem~\ref{th.54}, i.e., that~$\bO \in M$, follows by a simple argument 
from the results of Section~5, where we characterise all invariant sets that do not contain zero.  
The second part is more complicated. 
We give a proof of~Theorem~\ref{th.54} in subsection~6.3. Now we focus on sufficient conditions 
for convex eigensets. Let us recall that an invariant set~$M$ becomes  an eigenset if and only if 
every point~$\by \in M$ is reachable from~$M$ in arbitrary time (Definition~\ref{d.17}).  
 For convex sets,  
this condition can be relaxed to points~$\by \in \partial M$.
 \begin{prop}\label{p.13}
Suppose a convex body~$M$ contains zero and  is 
invariant for a normalized system~$\cA$; then, if every point of~$\partial M$ is reachable 
from~$M$ in arbitrary time, then~$M$ is an eigenset of the system~$\cA\cup \{-I\}$.   
\end{prop}
{\tt Proof}. 
We need to show that every point~$\bx \in {\rm int}\, M$
is reachable from~$M$ in arbitrary time. If~$\bx = \bO$, this is obvious: 
since~$\bO \in M$, this point is reachable from itself. 
If~$\bx \ne \bO$, then draw a ray~$\bO\bx$ that intersects~$\partial M$ at some point~$\bx_0$. 
The regime~$A(\vardot) \equiv -I$ generates a straight trajectory from~$\bx_0$  to~$\bx$. 
Since~$\bx_0$ is reachable in arbitrary time, so is~$\bx$ by Lemma~\ref{l.10}.   

{\hfill $\Box$}

 \begin{prop}\label{p.15}
Suppose a trajectory~$\bx(\vardot)$ is contained in a  convex body~$M$; 
then the dimension of the carrier of the point~$\bx(t)$ is non-decreasing in~$t$.   
\end{prop}
This means that if~$\by$ is reachable from~$\bx$, then~$\delta(\by) \ge \delta (\bx)$.  
 \smallskip 
 
{\tt Proof}. If a point~$\bx \in P$ is moved to a point~$\by$
by a switching law, then~$\by = \Pi \bx$, where~$\Pi$ is the corresponding 
fundamental matrix. Since~$\Pi$ is non-singular, the assertion follows.

{\hfill $\Box$}

\begin{cor}\label{c.38}
If~$M$ is a convex invariant body  of a normalized system~$\cA$,  
then a boundary point of~$M$ can be reached only from another boundary point.  
The corresponding trajectory lies on the boundary.   
\end{cor}
This corollary complements~Proposition~\ref{p.13}. 
The condition that every point from~$\partial M$
is reachable from~$M$ means that it is reachable from~$\partial M$
and all corresponding trajectories lie on the boundary.
The second  corollary of  Proposition~\ref{p.15} concerns extreme points 
of a convex eigenset. 
 \begin{cor}\label{c.40}
An extreme point~$\bx$ of a convex eigenset is reachable only from extreme points.   
 Every trajectory leading to~$\bx$ consists of extreme points. 
\end{cor}
We call an extreme point of a convex body {\em isolated} if it has a 
neighbourhood that does not contain other extreme points. 
 \begin{cor}\label{c.50}
An isolated extreme point of an eigenset~$M$ is reachable only from itself and 
belongs to a kernel of some 
operator from~$\cA$. 
\end{cor}
{\tt Proof}. If an extreme point~$\by$ is reachable from another 
extreme point~$\bx$ in some time~$t> 0$, then the corresponding trajectory
passes though any neighbourhood of~$\by$, which is impossible. 
Hence, $\by$ is reachable only from itself by a stationary trajectory
of some regime~$A \in \cA$. Hence,~$A\by  = 0$. 

{\hfill $\Box$}

\bigskip 	 

\begin{center}
\textbf{6.2. Polyhedral eigensets} 
\end{center}

\bigskip

If an eigenset~$M$ is a convex polytope, then we call it {\em eigenpolytope}. 
Corollaries \ref{c.38} -- \ref{c.50} being applied to an eigenpolytope 
mean the following:  every face~$L$ of an eigenpolytope, $0\le {\rm dim}\, L \, \le d$, 
is autonomous, i.e.,  all points from~$L$ are reachable from~$L$. It does not mean that 
$L$ is an eigenset since some trajectories starting on~$L$ may leave~$L$. 

 \begin{theorem}\label{th.54.5}
If a normalized  system~$\cA$ has an eigenpolytope~$M$, then every vertex of~$M$ belongs to a kernel of some operators from~$\cA$. 
Every  face~$L$ of~$M$ possesses the property:  
each point of~$L$ is reachable from~$L$ in arbitrary time. 
\end{theorem}
{\tt Proof}. The first assertion follows from Corollary~\ref{c.50}. 
The second one is proved by induction in the dimension~$s\, =\, {\rm dim}\, L$. 
For~$s=0$, i.e., when~$L$ is a vertex, the assertion hollows from Corollary~\ref{c.50}. 
Assume it holds for all faces of dimensions less than~$s$.
Every point from the boundary of~$L$ lies on a face of~$L$ of  
a smaller dimension, and hence, it is reachable from that face
(and hence, from~$L$) in arbitrary time. 
It remains to consider the case of an  interior point~$\by$ of~$L$. 
If it  is reachable 
from a point~$\bx\in M$ in time~$t$, then by Proposition~\ref{p.15}, 
$\delta (\bx) \le \delta (\by) = s$. Therefore, 
the trajectory leading from~$\bx$ to~$\by$ intersects the 
boundary of~$L$ at some point~$\bz$. 
By the argument above, $\bz$  is reachable from~$L$ in arbitrary time. 
Now Lemma~\ref{l.10} yields that~$\by$ is reachable from~$L$ in arbitrary time.
 
{\hfill $\Box$}

\bigskip


\bigskip 	 

\begin{center}
\textbf{6.3. Proof of Theorem~\ref{th.54}} 
\end{center}

\bigskip

We use the following simple consequence of Lemma~\ref{l.10}: 
 \begin{cor}\label{c.30}
 Let~$M$ be an invariant set of a normalized system~$\cA$
 and~$\bx_0 \in M$ be a point of a kernel of some operator from~$\cA$. 
 Then all points of all trajectories starting at~$\bx_0$
 are reachable from~$M$ in arbitrary time. 
 \end{cor}
{\tt Proof}. Let~$A\bx_0 = 0$
for some~$A\in \cA$; then the point~$\bx_0$ is reachable from~$\bx_0$ in arbitrary time 
by means of the regime~$A$. Applying Lemma~\ref{l.10} we complete the proof. 
 
{\hfill $\Box$}

\noindent {\tt Proof of Theorem~\ref{th.54}} is split into two parts. 

1. {\em Show that~$\bO \in M$.}  If this is not true, then by Theorem~\ref{th.50},   
~$M$ is a closure of the union of the sets~$rM^{(r)}$ over all~$r \in [0, +\infty)$, 
where each~$M^{(r)}$ is either empty or a compact subset of a unit sphere. 
The minimal~$r = r_{\min}$ for which~$M^{(r)} \ne \emptyset$ (it exists due to the compactness) is strictly positive. 
The intersection of~$M$ with the sphere of radius~$r_{\min}$ is an invariant set. 
If~$M$ is convex, then this intersection is a singleton. 
In this case, however, 
the straight line connecting  it with the origin 
is a common invariant subspace of the family~$\cA$, which is a contradiction. 
\smallskip 

2. {\em Prove that if~$\bO \in \partial M$, then~$\bO$ is a conical 
point}. Consider the maximal tangent subspace~$L$ at the point~$\bO$. 
This is a locus of points~$\bu$ such that the distance from the point~$\, \tau \, \bu$ to~$M$
is~$o(\tau)$ as~$\tau\to 0$. If~$\bO$ is not conical, 
then~${\rm dim}\, L \, \ge 1$. Denote by~$N$ the cone of outer normal vectors to~$M$ at~$\bO$
and~$\bar N = \{\bn \in N: \ \|\bn\| = 1\}$. Then~${\dim}\, N = d-k$ and~$L$ is an orthogonal complement 
to~$N$ in~$\re^d$. 
 Due to irreducibility of~$\cA$, there is a vector~$\bv\in L$
 and~$A\in \cA$ 
 such that~$A\bv \notin L$, i.e., for some~$\bn \in \bar N$, we have  $(\bn, A\bv) \ne 0$. 
 Denote~$(\bn, A\bv) = \alpha$. Taking, if necessary, 
 an opposite vector, we may assume that~$\alpha >   0$. 

For arbitrary~$\lambda > 0$, we introduce an operator~$R(\lambda)$
by the formula~$e^{\, \lambda A}\, = \, I\, + \, \lambda A \, + \, \lambda \, R(\lambda)$. 
Since  $R(\lambda) = o(\lambda)$ as~$\lambda \to 0$, one can choose a small $\lambda$ for which~$\|R(\lambda)\|\, <  \, 
\frac{\alpha}{2\, \|v\|}$. We fix this~$\lambda$. 
 
 For~$\tau \ge 0$, denote by~$\bb(\tau)$ the point of~$M$ closest to the point~$\tau\bv$. 
 We have~$\tau \bv \, - \, \bb(\tau) \, = \, o(\tau)$ as~$\tau \to 0$. Therefore, 
$$
e^{\, \lambda A} \bb(\tau) \ = \ 
\bb(\tau) \, + \, \lambda A\bb(\tau) \, + \, \lambda\, R(\lambda )\bb(\tau) \ = \ 
\tau \bv  \, + \, \tau \, \lambda A \bv \, + \, \tau \lambda\, R(\lambda )\bv\ \ + \ o(\tau)\quad 
\mbox{as}\ \tau \to 0\, . 
 $$ 
Taking into account that~$\bigl(\bn\, , \, \bv\bigr) \, = \, 0$, we obtain 
$$
\bigl(\, \bn\, , \, 
e^{\, \lambda A} \bb(\tau)\, \bigr) \ = \ 
\tau \, \lambda \, \bigl( \, \bn \, , \, A \bv\bigr)  \, + \,   \tau \lambda \bigl(\,\bn\, , \,  R(\lambda )\bv\, \bigr)
\ \ + \ o(\tau)\quad 
\mbox{as}\ \tau \to 0\, . 
 $$ 
Since, $(\bn\, , \, A\bv) \, =  \, \alpha$ and 
$\bigl(\bn\, , \, R(\lambda )\bv\, \bigr) \, \ge \, - \, \|\bn\|\, \|R(\lambda )\|\, \|\bv\|\, = \, 
- \,\|R(\lambda )\|\, \|\bv\| \, > \, - \, \frac{\alpha}{2}$, it follows that  
$$
\bigl(\, \bn\, , \, 
e^{\, \lambda A} \bb(\tau)\, \bigr) \ = \  \tau \lambda  \, \bigl[\, \alpha \, + \, 
 \bigl(\,\bn\, , \,  R(\lambda )\bv\, \bigr)\,  \bigr] \ + \ o(\tau) \quad \ge \quad 
 \frac{\tau \lambda \alpha}{2} \ + \ o(\tau) \qquad 
\mbox{as}\ \tau \to 0\, . 
$$
Thus, for small~$\tau > 0$, we have~$\bigl(\, \bn\, , \, 
e^{\, \lambda A} \bb(\tau)\, \bigr)\, > \, 0$, i.e., 
the point~$\bx(\lambda)$
of the trajectory~$\bx(t)\, = \, e^{t A}\bb(\tau)$ starting at~$\bb(\tau) \in M$ is outside of~$M$. 
This means that~$M$ is not invariant.

{\hfill $\Box$}

\bigskip

\begin{center}
\large{\textbf{7. The realization problem} }
\end{center}

\bigskip 

Having established the properties of eigensets, we arrive at a natural 
question: which compact set can be an eigenset? Given~$M \subset \re^d$, 
does there exist an irreducible  system~$\cA$ for which~$M$ is an eigenset? 
If the answer is affirmative, we say that {\em $M$ is realizable}, i.e., it can be realized as an eigenset
of a suitable irreducible system. 

We consider the realizability problem for convex sets. The answer depends on the geometry of~$M$ and on its location 
in~$\re^d$ (or, equivalently, on the  location of  the origin~$\bO$ with respect to~$M$). In the next section 
we focus on  the geometrical aspect only, regardless of the location. This is a different problem and, as we will see,  the answer can be nontrivial  already for flat polygons.

Theorem~\ref{th.54.5} implies that not every convex body~$M \subset \re^d$ is realizable. 
Indeed,   $M$ has to at least 
contain the origin, either as an interior or as a conical point. This condition  is, however, not sufficient. 
In this section,  the answer will be obtained   for some classes of convex sets.  

\bigskip 	 

\begin{center}
\textbf{7.1. Auxiliary results} 
\end{center}

\bigskip

Several minor results are collected here to provide 
tools for the construction of systems from given eigensets.
We begin with  the following technical lemma
that spares us from the need to  check the irreducibility of the system. We formulate it in a stronger form, 
for invariant set.

\begin{lemma}\label{l.30}
If a convex body~$M$ is invariant for a reducible system~$\cA$ and 
contains the origin either as an interior point or as a conical point, then 
~$\cA$ can be complemented by one operator to become irreducible and 
still have~$M$  invariant. 
\end{lemma}
{\tt Proof}. For an arbitrary element~$\bar A \in \cA$ which is not proportional 
to the identity, there is an operator~$C$ without common invariant subspaces with~$\bar A$. 
If~$\bO \in {\rm int}\, M$, then~$M$ is invariant for the operator~$\lambda C -  I$, 
whenever~$\lambda$ is small enough. Hence, $M$ is invariant for an irreducible 
system~$\cA \cup \{\lambda C -  I\}$. If~$\bO$ is a conical point, then 
the conic hull of~$M$ is a proper cone~$K$. There exists an operator~$C$ that 
maps~$K\setminus \{0\}$ to~${\rm int}\, K$.  After a small perturbation, this 
operator does not have common invariant subspaces with~$\bar A$. 
Then~$M$ is invariant for the operator~$\lambda C -  I$ (a {\em $K$-Metzler operator}, see~\cite{GLP17})
whenever~$\lambda > 0$ is small.

{\hfill $\Box$}

Lemma~\ref{l.30}  allows us not to think of irreducibility and even not to mentioned 
this property in the proofs. Note, however, that this is true only for convex eigensets which are, moreover,  
contain zero as an interior or as a conical point. Now we introduce a class of 
singular regimes which will be useful to construct systems with convex eigensets.

\begin{defi}\label{d.40}
The {\em compression to zero} is the operator~$A= - I$. 
For a given direct decomposition $\re^d = U \oplus V$, where~$U, V$ are subspaces 
of dimensions~$j$ and~$d-j$ respectively, $1\le j\le d-1$, 
the {\em compression towards~$U$ parallel to~$V$} is the operator~$A$
such that~$A|_{U} = 0, \, A|_{V} = -I$. 
\end{defi}
We identify an operator 
of compression with the corresponding regime and with the constant switching law defined 
by it. 
The compression to zero (or the corresponding switching law)
moves an arbitrary point~$\bx_0 \in \re^d$ to the origin along the straight trajectory~$\bx(t)\, = \, e^{-t}\bx_0, \, t\ge 0$. 
If we are given a direct sum $\re^d = U \oplus V$, then every point has a unique 
decomposition~$\bx = \bu + \bv$ with~$\bu \in 
U, \bv\in V$. The compression towards~$U$ parallel to~$V$ moves an arbitrary point~$\bx_0 = \bu_0 + \bv_0$
along the straight trajectory~$\bx(t) \, = \, \bu_0 \, + \, e^{-t}\bv_0$
to the limit point~$\bu_0$ which is a projection of~$\bx_0$ to the subspace~$U$ 
parallel to~$V$. This trajectory fills the half-interval~$[\bx_0, \bv_0)$

 \begin{prop}\label{p.18}
Let~$M$ be an arbitrary convex body that contains the origin; 
then~$M$ is invariant for the compression to zero. 

Let~$U$ be a straight line (one-dimensional subspace)
and~$V$ be a hyperplane not containing~$U$. Let~$U$ intersect~$M$ by a segment~$[\bu_1, \bu_2]$; then  
the following holds: ~$M$ is invariant for the compression towards~$U$ parallel to~$V$ 
 if and only if the  hyperplanes~$\bu_1+V$ and~$\bu_2 + V$ are both supporting for~$M$.  
\end{prop}
{\tt Proof}. The compression to zero moves every point~$\bx_0 \in M$
along the half-interval~$[\bx_0, \bO)$ which lies in~$M$. Hence, $M$ is invariant for this compression. 
The operator~$A$ of compression towards~$U$ parallel to~$V$ moves a point~$\bx_0 \in M$
along the half-interval~$[\bx_0, \bu_0)$, where~$\bu_0$ is a projection of~$\bx_0$ to the line~$U$
parallel to the hyperplane~$V$.  If the hyperplanes~$\bu_1+V$ and~$\bu_2 + V$ are both supporting for~$M$, 
then~$\bu_0 \in [\bu_1, \bu_2] \subset M$. Since~$M$ is convex,  the 
 whole segment~$[\bx_0, \bu_0)$ is 
contained in~$M$. Thus, the trajectory of the regime~$A$ starting at an arbitrary point~$\bx_0 \in M$
never leaves~$M$, consequently, $M$ is invariant for~$A$. 

Conversely, if, for example, the hyperplane~$\bu_1+V$ is not supporting, then 
there is a point~$\bx_0 \in M$ such that this hyperplane intersects the segment~$[\bO, \bx_0]$
at its interior point. This means that the projection~$\bu_0$ of the point~$\bx_0$
to the line~$U$ parallel to~$V$ lies outside~$M$. Hence, the trajectory of the regime~$A$
starting at~$\bx_0$ leaves the body~$M$. In this case $M$ is not invariant.

{\hfill $\Box$}

\bigskip 	 

\begin{center}
\textbf{7.2. Realizable convex bodies} 
\end{center}

\bigskip

We find several classes of realizable convex sets containing the  origin~$\cO$; among them are  
all simplices, parallelepipeds and Euclidean balls.  Some mild assumptions (not very obvious)
 must be satisfied. For example, a realizable polyhedron has to contain the 
 origin either inside or at a  vertex. All Euclidean balls are realizable in dimensions~$d\ge 3$, while 
 for $d=2$, realizable only those centered at the origin.

\begin{theorem}\label{th.55}
An arbitrary convex body symmetric about the origin is realizable.   
\end{theorem}
{\tt Proof}. Denote this  body by~$M$. 
For an arbitrary point~$\by \in \partial M$, 
draw an affine supporting hyperplane~$\by + V$ to~$M$, where~$V\subset \re^d$
is a linear hyperplane. By the symmetry, the hyperplane~$- \by + V$ is also supporting.
Invoking now Proposition~\ref{p.18} we obtain that~$M$ is invariant for  the compression~$A_{\by}$
towards the line~$U = \bO\by$ parallel to~$V$. For this compression, $\by$ is a fixed point and hence, 
is reachable from~$\by$ in arbitrary time. Since~$\by$ is an arbitrary point from~$\partial M$, we see that,  
for the family~$\cA = \{A_{\by}: \, \by \in \partial M\}$, 
all boundary points are reachable from~$\partial M$ in arbitrary time. 
By Proposition~\ref{p.13},  $M$ is an eigenset of the family~$\cA\cup \{-I\}$.

{\hfill $\Box$}

One may think that Theorem~\ref{th.55} implies, for examples, the 
realizability of balls and parallelepipeds, since they are symmetric. This is true, but only 
for the sets symmetric about the origin. Nevertheless, we are going to see that all parallelepipeds and balls 
containing the origin (non necessarily at the center) are indeed realizable. The same holds for simplices, which are already not symmetric.

\begin{theorem}\label{th.60}
An arbitrary simplex in~$\re^d, \, d\ge 2$, containing the origin  inside or at a vertex,  is realizable.  The same is true for a parallelepiped. 
\end{theorem}
{\tt Proof}. {\em Parallelepiped}. After a change of coordinates it can be assumed 
that~$M$ is a unit cube containing the origin and having its edges parallel to 
the coordinate axes. Thus,~$M = \{\bx \in \re^d: \ a_i \le x_i \le b_i, \, i=1, \ldots , d\}$, 
where~$b_i - a_i = 1, \, a_ib_i < 0$ for all~$i$.
Consider first the case~$\bO \in {\rm int}\, M$. 
A line~$\bO\by$ drawn through an arbitrary vertex~$\by$, intersects~$\partial M$
for the second time at a point~$\bz$. 
Denote by~$V$ the facet containing~$\bz$ (if there are several ones, take any of them). 
By Proposition~\ref{p.18}, the cube~$M$ is invariant for 
the compression towards the line~$\bO\by$ parallel to~$V$. For this operator, 
the vertex~$\by$ is fixed, hence, $\by$ is reachable  in arbitrary time.
Taking such operator for every vertex, we obtain a system for which 
all the vertices of~$M$ are reachable in arbitrary time.   

 To an arbitrary point~$\bx \in M, \, \bx \ne 0$, 
we associate a vertex~$\by$ of~$M$ as follows: 
if~$x_i = 0$, then~$y_i = 0$; if~$x_i < 0$, then~$y_i = a_i$; 
if~$x_i > 0$, then~$y_i = b_i$. Then 
consider a diagonal matrix~$A$ as follows: 
if~$y_i = 0$, then~$a_{ii} = 0$; otherwise,~$a_{ii} = \ln \, \frac{x_i}{y_i}$. 
Note that all elements of~$A$ are non-positive, hence,~$M$ is invariant for the regime~$A$
Then the trajectory
$\bx(\vardot)$ generated by~$A$ going from~$\bx(0) = \by$
is such that~$\bx(1) = \bx$. Thus,~$\bx$ belongs to the trajectory going from a point 
reachable from~$M$ in arbitrary time. By Lemma~\ref{l.10}, $\bx$ is also 
reachable in arbitrary time. The union of those regimes over all points~$\bx \in M\setminus \{\bO\}$
define the system for which~$M$ is an eigenset. 

It remains to consider the case when~$\bO$ is a  vertex of~$M$. 
Take arbitrary~$\bx \in M\setminus \{\bO\}$. Let the line~$\bO\bx$ intersect~$M$
by the segment~$[\bO, \bz]$ and~$V$ be a facet containing~$\bz$. Then we 
associate to~$\bx$ the operator of compression towards the line~$\bO\bx$ parallel to ~$V$. 
By this operator~$\bx$ is reachable in arbitrary time. Hence,~$M$ is an eigenset of the family of all those operators over all points~$\bx \in M\setminus \{\bO\}$.

\smallskip 

 {\em Simplex}. Consider first the case~$\bO \in {\rm int}\, M$. Induction in~$d$. 
 Let~$d\ge 3$ and the assertion is true in dimension~$d-1$. 
 For an arbitrary vertex~$\bv$ of~$M$, draw a line~$\bO\bv$ which intersects 
 the opposite facet~$V$ at a point~$\bv'$. Denote by~$L$ the affine hyperplane 
 containing~$V$ and make it a linear space by setting the origin at the point~$\bv'$. 
 Since~$V$ is a~$(d-1)$-dimensional simplex in~$L$, 
 by the inductive assumption, there is a system~$\cA_{\bv}$ of operators on~$L$
 for which~$V$ is an eigenset. To every operator~$A \in \cA_{\bv}$ we associate an operator~$\tilde A$ on~$\re^d$
 as follows:~$A\bv = 0$ and~$\tilde A|_{L} = A$. For the obtained system~$\tilde \cA_{\bv}$, 
 the set~$V$ is an eigenset and the set~$M$ is invariant. The union of those systems 
 over all vertices~$\bv$ of~$M$ is a system for which all boundary points of~$M$ are reachable 
 in arbitrary time. Hence, by  By Proposition~\ref{p.13}, $M$ is an eigenset. It remains to consider the base of induction~$d=2$. In this case~$M$ is a triangle, for which we do the same: 
For every vertex $\bv$, draw a line~$\bO\bv$ which intersects 
 the opposite facet~$V$ at a point~$\bv'$. Then consider  the operator of compression 
 towards the line~$\bO\bv$ along~$V$. Those three operators 
 form a system for which~$M$ is an eigenset. 
 
 The case when~$\bO$ is a vertex of~$M$ is considered in the same way as for the cube 
 (considered above).

{\hfill $\Box$} 
\medskip 

Thus, in all dimensions~$d\ge 2$,  each simplex and each parallelepiped is realizable provided that 
the origin~$\bO$ lies either in the interior or at a vertex. For other convex sets, 
the situation can be different. A $d$-dimensional ball 
containing the origin is always realizable, but only for~$d\ge 3$. Since the realizability property 
is invariant with respect to linear transforms, the same is true for ellipsoids. 
\begin{theorem}\label{th.65}
For~$d\ge 3$, an ellipsoid is realizable  if and only if it contains the origin as an interior 
point.  For~$d=2$, an ellipse is realizable if and only if it has center at the origin. 
\end{theorem}
{\tt Proof}. ($d\ge 3$) The necessity follows from Theorem~\ref{th.54} because an ellipsoid 
does no have conic points. 
To prove the sufficiency we assume that~$M$ is a unit Euclidean ball
and its center has coordinates $(c, 0, \ldots ,0) \in \re^d, \, c\in [0,1)$. 
 The operator~$A$ that maps a point~$(x_1, \ldots , x_d)$ to~$(0, x_3, - x_2, 0, \ldots , 0)$, defines 
 trajectories~$\bx(t) = e^{tA}\bx_0 \, = \, (x_1, y_2(t),  y_3(t),  x_4, \ldots , x_d), \, t\ge 0$, 
 where the vector~$(y_2(t),  y_3(t))$ is the image of the rotation of the vector~$(y_2(0), y_3(0)) = (x_2, x_3)$
 by the angle~$t$. Clearly, $M$ is the eigenset of the one-element system~$\cA = \{A\}$. 
\smallskip 

($d = 2$). The sufficiency follows from Theorem~\ref{th.55}, the necessity is proved in Proposition~\ref{p.50}
in Section~8.

{\hfill $\Box$} 
\medskip

What could be a criterion of realizability? We do not know the answer and leave an open problem: 
\begin{problem}\label{prob.10}
Characterize all realizable convex sets in~$\re^d$. Among them, find all convex
sets that are realizable for any location of the origin in their interiors. 
\end{problem} 
An example of a non-realizable convex set~$M$ such that~$\bO \in {\rm int}\, M$ is given by a disc in~$\re^2$
with center different from~$\bO$ (Theorem~\ref{th.65}). More examples including 
polygons are in the next section. Although every triangle and parallelogram 
is realizable by~Theorem~\ref{th.60}, there are non-realizable quadrilaterals. 
Moreover, every convex 
quadrilateral different from a parallelogram  can be  positioned on the plane to  
become non-realizable (although containing  
the origin inside). 

 \bigskip

\begin{center}
\large{\textbf{8. Plane eigensets}  }
\end{center}

\bigskip

We say that a compact set~$M$ {\em has a realizable copy }
 if there is a realizable set~$M'$ isometric to~$M$. 
Since the realizability is invariant with respect to linear transforms, 
we see that an isometry can be replaced merely by a parallel translation. 
Thus,~$M$  has a realizable copy if it becomes realizable after a translation. In other words, 
if there is a position of the origin~$\bO$ and an irreducible system~$\cA$
for which~$M$ is an eigenset.  

By Theorem~\ref{th.60}, for  a triangle or a parallelogram, every 
copy containing zero inside is realizable. However, for a disc, this is not the case. 
Not every quadrilateral is realizable but it always possesses 
a realizable copy, this is shown in  Proposition~\ref{p.50} below. 
The main result of this section,  Theorem~\ref{th.70}, states that 
every $n$-gon for~$n\le 5$, has a realizable copy, while 
for hexagons, this is not true. We begin with a preliminary work.

 \newpage 

\begin{center}
\textbf{8.1. Auxiliary facts and notation}  
\end{center}

\bigskip

A system is {\em non-defective} if there is a constant~$C$ such that 
for every trajectory, we have~$\|\bx(t)\|\, \le \, C\, \|\bx(0)\|, \, t\ge 0$.
It is known that a system is non-defective if and only if it possesses 
a convex  invariant body~\cite{MP86}.  
An irreducible system is always non-defective~\cite{B88}.   A matrix~$A$ is non-defective if 
every its eigenvalue~$\lambda$ such that~${\rm Re}\, \lambda \, = \, \sigma(A)$ 
is non-defective, i.e., possess only trivial Jordan blocks. 

We call a matrix~$A \in \cA$ {\em dominant} if~$\sigma(A)\, = \, \sigma(\cA)$. 
The stationary switching law~$A$ corresponds to the fastest growth of trajectories. 
All dominant matrices  of a non-defective systems are non-defective either. 
In particular, for all singular matrices of an irreducible normalized 
system~$\cA$, zero is a non-defective eigenvalue. 
For~$d=2$, this means that every singular matrix~$A \in \cA$ 
 has the second eigenvalue negative. Indeed, it is real,  
cannot be positive  (since the system is normalized) and cannot be zero
because~$A\ne 0$. Thus, we have proved
  \begin{prop}\label{p.48}
Every singular matrix~$A$ of an irreducible normalized 2D system defines a compression 
towards~${\rm Ker}\, A$ along~${\rm Im}\, A$.    
\end{prop}

 \begin{lemma}\label{l.20}
Let  a normalized system~$\cA$ have a convex invariant set~$M$ and some~$A\in \cA$ is singular. 
Let the kernel of~$A$ intersect~$M$ by a segment~$[\ba, \bb]$; 
then the two lines drawn through~$\ba$ and~$\bb$  parallel to~${\rm Im}\, A$
are supporting for~$M$. 
\end{lemma}
Thus, a convex invariant set lies in the strip parallel to~$\, {\rm Im}\, A$ passing through the points of intersection of~$M$ 
with~$\, {\rm Ker}\, A$. This is true for all singular operators from~$\cA$. 
\smallskip 

 {\tt Proof of Lemma~\ref{l.20}}. If the line~$\ell$ passing though~$\ba$ parallel to~$\, {\rm Im}\, A$ is not supporting for~$M$, then it 
strictly 
separates some point~$\bx_0 \in M$ from the origin. The trajectory of the 
regime~$A$ starting from~$\bx_0$ 
goes along~$\ell$ towards the line~$\ba\bb$ to their common point. 
Since this common point does not belong to~$M$, it follows that 
the trajectory leaves~$M$, which is impossible. 
 
{\hfill $\Box$}

 \bigskip 	 

\begin{center}
\textbf{8.2. Non-realizable 2D convex sets}  
\end{center}

\bigskip

  \begin{prop}\label{p.50}
A Euclidean disc centered not at the origin 
cannot be an invariant set of an irreducible system. 
\end{prop}
{\tt Proof}. Assume that a unit disc~$M$ with the center at a point~$\bc = (c,0), \, c > 0$, 
is invariant for~$\cA$. Denote by~$\ba\bb$ the diameter of~$M$
on the axis~$\bO X$. By Theorem~\ref{th.54}, the origin~$\bO$ is inside~$M$, 
hence,~$\bO$ belongs to the open interval~$(\ba, \bb)$. If~$\cA$ contains a singular 
operator~$A$, then its kernel intersects~$M$ by a chord and the lines passing through 
its ends  parallel to the image of~$A$ must be tangents of~$M$. 
We see that the tangents at the ends of a chord are parallel, therefore, 
this chord is a diameter passing though the origin, i.e., 
 it coincides with the segment~$\ba\bb$. Thus, {\em every singular operator}
 from~$\cA$ (if any) has the kernel~$\bO X$ and the image~$\bO Y$.  
 By~\cite[Theorem~2]{PM25}, if for a normalized system~$\cA$, no points outside~$\bO X$ belong to the kernels of operators from~$\cA$, 
 then there is a trajectory~$\bx(t)$ from an arbitrary point~$\bx_0 \in \bO X$  to a symmetric point~$-\bx_0$.
 Taking~$\bx_0 = \bb$, we obtain a trajectory to~$-\bb$. This is impossible, since 
 $-\bb \notin M$. 
 
{\hfill $\Box$}

\begin{prop}\label{p.45}
\textbf{a)} A  quadrilateral with diagonals intersecting at the origin
is realizable; 

\textbf{b)} Every convex 
quadrilateral, other than a parallelogram,  can be positioned in the plane so that it contains 
the origin inside but is not an eigenset 
of any system. 
\end{prop}
{\tt Proof}. Let a quadrilateral~$M = \ba\bb\bc\bd$ has its diagonals intersecting at~$\bO$.
The compressions towards~$\bO\ba$ and~$\bO\bb$ along~$\bO\bb$ and~$\bO\ba$
respectively, leave the set~$M$ invariant (Proposition~\ref{p.18}). 
This implies that all the vertices  
are reachable from themselves in arbitrary time. 

Take now the side~$\ba\bb$ and choose the most distant vertex 
from the line containing this side. Let it be~$\bc$. Denote~$\bv = \bb - \ba$. 
The compression towards~$\bO\ba$ along~$\bv$  leaves~$M$ invariant.  
All points of the side~$\ba\bb$  are reachable by this regime in arbitrary time. 
Applying the same argument for other sides, we obtain that all points of the 
boundary of the quadrilateral are reachable in arbitrary time. 
By Proposition~\ref{p.13},  $M$ is an eigenset. 

Now prove that, unless the quadrangle $M$ is a parallelogramm, it can be situated on the 
plane so that 1) it has the origin in its interior;  2) there is no system for which it is an eigenset.
Take non-parallel opposite sides, let them be~$\ba\bb$ and~$\bc\bd$. 
Assume that~$\bc$ is more distant from the line~$\ba\bb$ than~$\bd$. 
Place the origin~$\bO$ in an arbitrary interior point of the triangle~$\ba\bb\bd$
and prove that~$M$ is not an eigenset of any system. 
If such a system~$\cA$ exists, then by~Theorem~\ref{th.54.5}, the line~$\bO\bd$ is a kernel of some operator~$A\in \cA$. 
Denote by~$\bv$ the point of intersection of this line with the side~$\ba\bb$. 
Then by Lemma~\ref{l.20}, there is a pair of parallel lines of support of~$M$ passing through 
$\bv$ and~$\bd$. There is a unique supporting line at~$\bv$, this is~$\ba\bb$. 
Hence, the line parallel to~$\ba\bb$ passing through~$\bd$ is also supporting. 
This is impossible since it intersects the interior of~$M$.

{\hfill $\Box$}

\bigskip 	 

\begin{center}
\textbf{8.3 Realization of plane polygons. The fundamental theorem}  
\end{center}

\bigskip 

Every triangle, quadrangle, and pentagon has an isometric  copy which is an eigenset.  On the other hand, there are hexagon  which do not have isometrical eigensets. 
\begin{theorem}\label{th.70}
For~$n=3,4,5$, every convex~$n$-gone on the plane possesses a realizable copy.   
For~$n=6$, this is not true:    
there is a hexagon that does not have a realizable copy. 
\end{theorem}
{\tt Proof}. The cases~$n=3$ and~$n=4$ follow from~Theorem~\ref{th.60} and Proposition~\ref{p.45} respectively. 

$\mathbf{n=5}.$ Firstly we define some notation. Each side of a convex pentagon~$\ba\bb\bc\bd\be$ has a unique opposite 
vertex, which is not a vertex of the neighbouring sides. 
A vertex is called {\em normal} if it is  the most distant 
vertex from the line containing the opposite side. The inequality is not strict: 
one of the two neighbouring vertices may have the same distance (in which case 
the two corresponding sides of the pentagon are parallel). The side opposite to the 
normal vertex is also called normal. A side~$\bb\bc$ is called {\em short} 
if the sum of the interior angles at the vertices~$\bb$ and~$\bc$ is 
at least~$\pi$. This is equivalent to say that the lines of the 
two neighbouring sides are either parallel or intersect 
on the other side with respect to~$\bb\bc$ than the pentagon.   
If the side~$\ba\bb$ is normal, then the 
ray  $\bd\bc$ either intersects the ray~$\ba\bb$ of is parallel to it. 
Hence, the side $\bb\bc$ is short. Thus, {\em both of the neighbouring  sides of a normal side are short}. 
The converse is also true: the side between two short sides is normal. 

We are constructing a system~$\cA$ that 
has a given pentagon~$M \, = \, \ba\bb\bc\bd\be$ as an eigenset.
Two cases need to be considered: 
\smallskip 

\begin{figure}[h!]
\begin{center}
\includegraphics[width=0.35\linewidth]{Illustrations1.mps}
 	\caption{{\footnotesize Pentagon, case 1}}
 	\label{f.10}
 \end{center}
 \end{figure}

1. There are three successive normal sides. In this case, all five sides are normal. 
  Indeed, if the sides~$\ba\bb, \bb\bc, \bc\bd$ are normal, then all 
  their neighbouring sides, i.e., all sides of~$M$ are short, hence, all 
  sides are normal. Choose the point~$\bO$ arbitrarily from the 
  convex pentagon formed by  the diagonals of~$M$, Fig.~\ref{f.10}. Then the line~$\bd\bO$  intersects 
  the opposite side~$\ba\bb$ at a point~$\bd'$. Let~$A$ be the 
  compression towards~$\bd\bd'$ along~$\ba\bb$. 
  It leaves the set~$M$ invariant and leaves the segment~$[\bd, \bd']$ 
  fixed. Let~$\cA'$ be the system of five such compressions (for all sides of~$M$). 
  Each vertex is a fixed point for a suitable compression and hence, it is reachable 
  in arbitrary time. In particular, $\ba, \bb$ are both reachable 
  in arbitrary time. Hence, so are all points of the side~$\ba\bb$ (by means of the 
  compression~$A$). We see that all points of the boundary of~$M$ are reachable in arbitrary time, therefore, by Proposition~\ref{p.13},~$M$ is an eigenset of the family~$\cA = \cA'\cup \{-I\}$.  
  \smallskip

2. The alternative for the case 1 is that there are two non-neighbouring sides wich are both not normal. Let them be~$\bb\bc$ and~$\be\ba$. Now we have to consider two subcases: 
 \smallskip

2a. The side~$\ba\bb$ is short, Fig.\ref{fig20a} In this case we set~$\bO = \bd$ and introduce the  
following system~$\cA$. The operator~$A_1$ 
defines a compression towards~$\bO\be$ along~$\ba\be$. Similarly,~$A_2$ is 
the compression towards~$\bO\bc$ along~$\bb\bc$. For those operators, the pentagon $M$ is invariant. 
Indeed, for both sides~$\ba\be$ and~$\bb\bc$, the most distant vertex is~$\bd = \bO$, 
hence, the lines drawn through~$\bO$ parallel to those sides are supporting for~$M$. 
Define also~$A_3, A_4$ as compressions along~$\ba\bb$ towards~$\bO\ba$ and~$\bO\bb$  respectively. 
Since the line drawn through~$\bO$ parallel to~$\ba\bb$ is supporting for~$\bO$, 
it follows that~$M$ is invariant for~$A_3$ and~$A_4$. All points of the boundary are reachable 
in arbitrary time, hence $M$ is an eigenset of~$\cA = \{A_1, \ldots , A_4, -I\}$. 
   \smallskip
   
   \begin{figure}[htbp]
  \centering
  \begin{minipage}[b]{0.35\textwidth}
    \includegraphics[width=\linewidth]{Illustrations2a.mps}
    \caption{Pentagon, the case 2a}
    \label{fig20a}
  \end{minipage}
  \hfill
  \begin{minipage}[b]{0.35\textwidth}
    \includegraphics[width=\linewidth]{Illustrations2b.mps}
    \caption{Pentagon, the case 2a}
    \label{fig20b}
  \end{minipage}
\end{figure}
   
2b. The side~$\ba\bb$ is not short, i.e., the rays~$\ba\be$ and~$\bb\bc$ intersect at some point~$\bu$. 
In this case both of  the rays $\bc \bd$ and~$\be\bd$ intersect the opposite sides of the 
triangle~$\bc\bu\bd$, Fig.\ref{fig20b} This implies that both of  the corresponding sides  $\bc \bd$ and~$\be\bd$ of~$M$ are short. 
If we assume that the side~$\bc\bd$ is not normal, 
then we return to the case 2a. Indeed, we have two not normal sides ($\bc\bd$ and~$\be\ba$)
separated by a short side ($\bd\be$). The same is when~$\bd\be$ is not normal. 
Thus, it remains to consider the case when both~$\bc\bd$ and~$\bd\be$ are normal. 
We set~$\bO\, = \, \ba\bc \, \bigcap \, \bb\be$ and take the following five operators: 
compressions towards~$\bb\be$ along~$\ba\be$ and~$\bd\be$; 
compressions towards~$\ba\bc$ along~$\bb\bc$ and~$\bd\bc$; 
the compression towards~$\bO\bd$ along~$\ba\bb$. 
The set~$M$ is invariant for all those operators and every point of the boundary
is reachable in arbitrary time,   hence,~$M$ is an eigenset for~$\cA = \{A_1, \ldots , A_5, -I\}$.

\medskip 

$\mathbf{n=6}.$  We give an example of a hexagon which cannot be presented as 
an eigenset of an irreducible system, even if we can choose the location of the origin. 
The following elementary fact will be used:  
\smallskip 

{\em If for a side~$\bu\bv$ of a convex polygon~$P$, 
the sum of interior angles at the vertices~$\bv, \bu$ is greater than~$\pi$, then 
there is no pair of parallel supporting lines of~$P$ passing through those vertices. }
\smallskip 

\begin{figure}[h!]
\begin{center}
\includegraphics[width=0.4\linewidth]{Illustrations3.mps}
 	\caption{{\footnotesize The hexagon which cannot be realized as an eigenset}}
 	\label{f.30}
 \end{center}
 \end{figure}

The same is true if~$\bu\bv$ is a diagonal of~$P$. {\em If the sum of two angles 
made by~$\bu\bv$ with two sides going from~$\bu$ and~$\bv$ in one half-plane with respect to~$\bu\bv$
is greater than~$\pi$, then 
there is no pair of parallel supporting lines of~$P$ passing through~$\bu$ and~$\bv$.}
\smallskip 

Consider an isosceles trapezoid~$\ba\bb\bc\bff$ with the biggest base~$\ba\bb$. 
Construct externally on the side~$\bff\bc$ an isosceles trapezoid~$\bff\bc\bd\be$ 
with the  biggest base~$\bff\bc$ and such that the angle~$\angle \, \bff\bc\bd$
is smaller than~$\angle \, \ba\bb\bc$. Finally, slightly 
move the vertex~$\bd$ towards the line~$\bff\bc$, Fig.\ref{f.30}.  The obtained convex 
hexagon~$M=\ba\bb\bc\bd\be\bff$ possesses the desired property. 
If, to the contrary,~$M$ is an eigenset of a system~$\cA$, then, by Theorem~\ref{th.54},  
the origin~$\bO$ can be either an interior point or a vertex of~$M$. 
Each of the three sides~$\bb\bc$, $\bd\be$, and~$\bff\ba$ possesses the property that  the 
sum of two adjacent angles of the hexagon is greater than~$\pi$. If~$\bO$
coincides with a vertex, say,~$\ba$, then the corresponding side~$\ba\be$
must be a kernel of some operator~$A\in \cA$ (Theorem~\ref{th.54.5}). 
Hence, by~Lemma~\ref{l.20}, there exist two parallel supporting lines of~$M$
passing through the ends of the side~$\ba\bb$, which is impossible. Thus,~$\bO$ cannot be a vertex. 

If~$\bO \in {\rm int}\, M$, then the line~$\bff\bO$  intersects~$M$ by a segment~$\bff\bff'$, 
where~$\bff'$ lies on one of the sides~$\ba\bb, \bb\bc, \bc\bd, \bd\be$. 
 Lemma~\ref{l.20} yields that there exist two 
 parallel lines of support for~$M$ passing through~$\bff$ and~$\bff'$. 
If~$\bff'$ is an interior point of a side of~$M$, then 
both those lines of support are parallel to this side.  
However, all lines passing through~$\bff$ parallel to the four aforementioned sides 
intersect the interior of~$M$ and hence, cannot be supporting. 
This proves that~$\bff'$ must be a vertex, i.e., one of the 
vertices~$\bb, \bc, \bd$. The case~$\bff' = \bc$ is impossible, since 
$\angle \ba\bff\bc + \angle \bb\bc\bff \, > \, \pi$ and therefore, 
there are no parallel supporting lines passing through~$\bff$ and~$\bc$. 
The case~$\bff' = \bd$ is impossible either, since 
$\angle \ba\bff\bd + \angle \bc\bd\bff \, > \, \pi$ and 
there are no parallel supporting lines passing through~$\bff$ and~$\bd$.
Therefore,~$\bff' = \bb$. Thus, $\bO \in \bff\bb$. 
Applying the same argument to the vertex~$\bc$ instead of~$\bff$, we show that 
 $\bO \in \bc\ba$. Consequently,~$\bO$ is a point of intersection 
of the diagonals~$\bff\bb$ and~$\bc\ba$. In this case the line~$\bd\bO$
intersects the side~$\ba\bb$ in its interior point, 
therefore, the line drawn through the vertex~$\bd$ parallel to the side~$\ba\bb$
must be supporting for~$M$ (Lemma~\ref{l.20}). This is not true, 
since this line strictly separates~$\be$ from~$\ba\bb$ and hence, intersects the interior of~$M$. The contradiction completes the proof.

{\hfill $\Box$}

 \bigskip 	 

\begin{center}
\large{\textbf{9. The problems of uniqueness and of convexity}  }
\end{center}

\bigskip 

Already for planar systems ($d=2$),  there may be different situations 
with uniqueness and convexity of eigensets. There is  
 a large class of planar systems having a unique (up to multiplication by a constant) eigenset, which is, moreover, convex.  
 On the other hand, there are systems that have no convex eigensets
 but have infinitely many different non-convex ones.

Recall that a matrix~$A\in \cA$ is  dominant if~$\sigma(A) = \sigma(\cA)$. 
The following theorem was established in~\cite{PM25}: 
\smallskip 

\noindent \textbf{Theorem~B}. {\em If a system~$\cA$ in~$\re^2$ does not have a dominant regime, 
then it possesses a unique up to multiplication invariant norm. The unit sphere of that norm 
is a closed periodic trajectory of the normalized system~$\cA-\sigma(\cA)I$. }
\smallskip 

Thus, if a second order system does not have a dominant regime, then its invariant ball is 
defined by a unique extremal trajectory. This ball turns out to be a unique eigenset of this system. 
Let us note that, in general,  the invariant ball of a system is  not its eigenset  (Remark~\ref{r.20}). 
\begin{theorem}\label{th.80}
If a system in~$\re^2$ does not have a dominant regime, then  its unique, up to homothety, 
eigenset is the unit ball of its invariant norm. 
\end{theorem}
{\tt Proof}. Let the system be normalized,~$f$ be its invariant norm,~$\bar A(\vardot)$ 
be the periodic switching law for its trajectory that is  a unit sphere of~$f$. 
Consider an arbitrary eigenset~$M$ of this system and take a point~$\bx \in M$ with the maximal value of~$f$. 
Without loss of generality it can be assumed that~$f(\bx) = 1$. The trajectory with the 
periodic switching law~$\bar A(\vardot)$ starting at~$\bx$ describes the invariant sphere~$S$. 
Thus,  $S \subset M$ and~$M$ lies inside~$S$. 
Therefore~$S = \partial M$. Since the system does not possess a dominant matrix, 
it follows that for every~$A_0 \in \cA$, we have~$\sigma(A_0) < 0$. 
Therefore, the trajectory of the switching law~$A_0$ starting at~$\bx$
approaches to zero. Applying the periodic switching law~$\bar A(\vardot)$  to each point of that trajectory 
we obtain a set of invariant spheres located inside~$M$. Those spheres 
form a unit ball of the invariant norm, which, therefore, coincides with~$M$.

{\hfill $\Box$} 

The following example demonstrates the opposite situation
when a planar system does not have convex eigensets but has infinitely many non-convex ones. 
This means, in particular, that a convex hull of an eigenset is not necessarily an eigenset. 

\begin{ex}\label{ex.20} (A system all of  whose  eigensets are non-convex). 
{\em Let us show that the following system~$\cA = \{A_1, A_2\}$ of two  $2\times 2$ matrices does not possess  convex eigensets although has infinitely many non-convex ones (non homothetic to each other): 
$$
A_1 \ =  \ 
\left( 
\begin{array}{rc}
-1 & p\\
0 & 0
\end{array}
\right)\, , \qquad 
A_2 \ =  \ 
\left( 
\begin{array}{rr}
0 & 0\\  
-1 & -p
\end{array}
\right)\, ,   
$$
where~$p\in (0,1)$. A kernel of~~$A_1$ is the line~$r_1:  \  x_1 = px_2$, 
its  image is the axis~$\bO x_1$; the kernel of~$A_2$ is the line~$r_2:  \  x_2 = - px_1$, 
its image is the axis~$\bO x_2$. Hence, $A_i$ defines the compression towards~$r_i$ along~$\bO x_i, \, i=1,2$. 

Suppose~$\cA$ has a convex eigenset~$M$; then by~Theorem~\ref{th.54}, $\bO \in M$. 
Denote by~$\bp_i\bq_i$ the segment of intersection of~$M$ with~$r_i, \, i = 1,2$. 
Proposition~\ref{p.18} yields that the two horizontal lines passing through~$\bp_1, \bq_1$
are both supporting for~$M$, the same with vertical lines through~$\bp_2, \bq_2$. 
Thus,~$M$ is inscribed in the rectangle~$\bb\be\bd\bc$, Fig.\ref{f.40}. 
Since~$A_1,A_2$ are compressions towards the segments~$\bp_1\bq_1$ and $\bp_1\bq_1$ respectively,  
it follows that~$M$ contains those segments.

\begin{figure}[h!]
\begin{center}
\includegraphics[width=0.4\linewidth]{Illustrations4.mps}
 	\caption{{\footnotesize Construction of eigensets of the system~$\{A_1, A_2\}$}}
 	\label{f.40}
 \end{center}
 \end{figure}

Let~$\bh$ denote the point of intersection 
of the segment~$O\bd$ with the boundary of~$\partial M$. The supporting line~$\ell$ 
of the set~$M$ drawn through~$\bh$
intersects the segments~$\bd\bq_1$ and $\bd\bq_2$. Observe that 
both of the vectors~$A_i\bh, \, i=1,2$, starting at~$\bh$ are directed from the 
line~$\ell$, i.e., increase the distance to~$\ell$. Hence, this holds 
for all points of~$M$ close to~$\bh$. Therefore,~$\bh$ is unreachable from other points of~$M$ which 
is a contradiction. 

It remains to show that~$\cA$ possesses infinitely many eigensets. 
Choose arbitrary segments~$p_i\bq_i \subset r_i, \, i = 1,2$, intersected at zero 
and define~$M$ as the closure of the set reachable from them by the system~$\cA$. 
This way we can construct infinitely many non-homothetic eigensets. 
}
\end{ex}

\begin{remark}\label{r.30}
{\em Theorem~\ref{th.80} cannot be applied to the system constructed in Example~\ref{ex.20} 
because it has two dominant regimes: $A_1$ and $A_2$. The eigenset consists of four ivies with one-point 
stems~$\{\bp_i\}, \{\bq_i\}, \ i=1,2$ (see Theorem~\ref{th.30}).  
}
\end{remark}

\begin{center}
\large{\textbf{Appendix}}
\end{center}

\begin{remark}\label{r.100}
{\em If~$\cA$ is convex, then all reachable sets are compact. Indeed, the boundedness is obvious and 
the closeness follows from the existence of solution of the optimal control problem:  
\begin{equation}\label{eq.opt}
\left\{
\begin{array}{l}
f_0(\bx(t)) \, = \, \|\bx(t) - \by\|^2 \  \to \ \min; \\
\dot \bx(\tau ) \, = \, \varphi(\tau, \bx, A)\, = \, A(\tau) \bx(\tau),\\
A(\tau) \in \cA,\ \tau \in [0,t], \\ 
 \bx(0) \in M\, . 
\end{array}
\right. 
\end{equation} 
Here the function~$\varphi$
is linear both in~$\bx$ and in the control variable~$A$, 
which takes values from a convex compact set;   
the objective function~$f_0$ is continuous and convex. 
Under these conditions, the solution exists (see, for instance,~\cite{F62}). 
Thus, if there exists a sequence of trajectories~$ \bx_k$ such that the values~$\bx_k(t)$ 
converge to some point~$\by$, then there exists a trajectory for which~$\bx(t) = \by$. 
We see that~$M_t$ is closed.} 
\end{remark}

{\tt Proof of Theorem~\ref{th.20}}. Firstly, we consider  the case of an irreducible system. 
After a suitable translation, it can be assumed that
 $\sigma(\cA) = 0$. Then, by Corollary~\ref{c.10}, for every eigenset~$M$, we have  $M_t = M, \  t \in \re_+$. 
 
Consider the invariant ball~$B$ and the sphere~$S = \partial B$ of the system~$\cA$.
For an arbitrary point~$\bx_0 \in S$, the notation~$P_t$ denotes the closure of the set of points   
reachable from~$\bx_0$ in time at least~$t$. Clearly~$P_t \subset B$  for every~$t\ge 0$. 
Besides,~$P_t$ contains all points of the ``tail'' of extremal trajectory~$\{\bar \bx (\tau): \, \tau \ge t\}$
located on~$S$, therefore, $P_t \cap S \ne \emptyset$. By the nested compact set theorem,  the intersection~$M = \bigcap\limits_{t\ge 0}P_t$ is not empty. It, moreover, intersects the invariant sphere~$S$ 
with the set~$\bigcap\limits_{t\ge 0}(P_t\cap S)$, 
hence,~$M\ne \{0\}$. Let us show that~$M$ is an emptyset.  This set consists of points~$\by$
that are the limits of sequences~$\by_k(t_k)$ as~$k\to \infty$, where  
$t_k \to +\infty$ and $\by_k(\vardot)$ is a trajectory going from the point~$\by_k(0) = \bx_0$. 
For every~$h > 0$, the set~$M_h$ consists of points~$\Pi(h) \by, \, \by \in M$, which are   
limits of the sequences~$\Pi(h)\by_k(t_k)\, = \, \by_k(t_k+h)$ as~$k\to \infty$.
Here~$\Pi(\vardot)$ is a fundamental matrix of an arbitrary switching law. 
Since~$t_k + h \to \infty$ as~$k\to \infty$, it follows that~$\Pi(h)\by \in M$. 
Therefore,
$M_{h} \subset M$. For the inverse embedding, we pass to a subsequence and assume that~$t_k > h$ 
for all~$k$ and  that~$\by_k(t_k-h)$ converge to some point~$\bz \in M$.
We see that there exists a sequence of trajectories~$\bx_k$ emanating  from  
the point~$\bz$ for which~$\bx_k(h) \to \by$ as~$k\to \infty$. 
Relying  again on the existence of the solution for problem~(\ref{eq.opt}), 
we conclude that there is a  trajectory such that~$\bx(0) = \bz, 
\bx(h) = \by$.  Thus, every point~$\by \in M$ is reachable in time~$h$ from some point~$\bz \in M$. 
This yields that $M_h=M$, which completes the proof in the irreducible case. 

In the case of reducible~$\cA$ we apply the induction in the dimension~$d$. 
Let~$L$  be the minimal by inclusion 
common nontrivial invariant subspace of matrices from~$\cA$. 
In a suitable basis those matrices get a block upper triangular form:
$$
A \ =  \ 
\left( 
\begin{array}{ll}
A^{(1)} & C\\
0 & A^{(2)}
\end{array}
\right)\, , 
$$  
where the block~$A^{(1)}$ corresponds to the subspace~$L$ and~$A^{(2)}$ corresponds to some 
 subspace~$K$ that complements~$L$ to the whole~$\re^d$. 
Let~$\cA^{(1)}, \cA^{(2)}$ denote the corresponding matrix families. 
By the minimality of~$L$, it follows that~$\cA^{(1)}$ is irreducible, 
hence, it possesses at least one eigenset~$M^{(1)} \subset L$
and for all its eigensets, the eigenvalue is equal to~$\sigma(\cA^{(1)})$. 
All eigensets of~$\cA^{(1)}$  are also eigensets of~$\cA$, hence, 
$\sigma(\cA^{(1)}) = \alpha_j$ for some~$j \le r$. 
It is known that 
for block upper triangular matrices,  $\sigma(\cA)\, = \, \max \, \bigl\{\sigma(\cA^{(1)}), \sigma(\cA^{(2)}) \bigr\}$.
Each trajectory~$\bx(t)$ of the system~$\cA$
is presented in the form~$\bx(t)\, = \, \by(t)+\bz(t)$, where~$\by(t) \in L, \, \bz(t) \in K$
for all~$t$ and 
\begin{equation}\label{eq.red}
\left\{ 
\begin{array}{ccl}
\dot \by  &  =  & A^{(1)}(t)\, \by \, + \, C(t)\bz\\
\dot \bz  &  =  & A^{(2)}(t)\, \bz
\end{array}
\right.   
\end{equation}
This implies that if~$M$ is an eigenset of~$\cA$, then 
$\bar M \, = \, \{\bz:  \ (\by, \bz) \in M\}$ is an eigenset of~$\cA^{(2)}$ 
with the same eigenvalue. Thus, all eigenvalues of~$\cA$ are also 
eigenvalues of~$\cA^{(2)}$. By the inductive hypothesis, 
the latter has at most~${\rm dim}\, K \, \le \, d-1$ eigenvalues, the largest of which is~$\sigma(\cA^{(2)})$. 
Hence,~$\cA$ has at most~$d$ eigenvalues, 
the largest of which does not exceed~$\max\, \{\sigma(\cA^{(1)}), \sigma(\cA^{(2)})\} \, = \, 
\sigma(\cA)$. It remains to show that~$\sigma(\cA)$ is an eigenvalue of~$\cA$, i.e., 
it is realized in one of its eigensets. If~$\sigma(\cA^{(1)}) \ge  \sigma(\cA^{(2)})$, then this is obviously true:~$M^{(1)}$ is an eigenset of~$\cA$ with the  eigenvalue~$\sigma(\cA^{(1)}) =  \sigma(\cA)$. 
If~$\sigma(\cA^{(1)}) < \sigma(\cA)$, then~$\sigma(\cA^{(2)}) = \sigma(\cA)$ and by the induction hypothesis~$\cA^{(2)}$ has an eigenset~$\bar M \subset K$ with the eigenvalue~$\alpha = \sigma(\cA^{(2)}) = \sigma(\cA)$. 
It may not  be an eigenset for the system~$\cA$ since 
$K$ is not an invariant subspace of~$\cA$. After a normalization, we assume that
 $\sigma(\cA) = 0$ and so~$\sigma(\cA^{(2)}) = 0, \,  \sigma(\cA^{(1)}) < 0$.

Since~$\sigma(\cA^{(1)}) < 0$, we see that there are constants~$C_1, C_2$
such that for an arbitrary switching law, the following holds for equation~(\ref{eq.red}):  
if~$\|\by(0)\| \le C_1$ and~$\|\bz(\vardot)\| \le C_1$, then
$\|\by(\vardot)\| \le C_2$. This yields that all trajectories of the system~$\cA$
emanating from points of the set~$\bar M$ are bounded uniformly. 
Let~$P_t$ denote the closure of the set of all points reachable from~$\bar M$ 
in time at least~$t$. It follows by applying the theorem of nested compact sets that the intersection~$M = \bigcap\limits_{t\ge 0}P_t$
is not empty. It  cannot coincide with zero either because the projection of this set 
 to the subspace~$K$ parallel to~$L$ coincides with~$\bar M$. 
This implies that~$M$ is an eigenset of~$\cA$. The proof is literally the same as in the irreducible case.

{\hfill $\Box$}

 \end{document}